\documentclass[12pt,fleqn]{article}

\usepackage{amsmath}
\usepackage{amsthm}
\usepackage{amssymb}
\usepackage{amsfonts}
\usepackage{graphicx}
\usepackage{hyperref}
\usepackage{appendix}

\newtheorem{lemma}{Lemma}

\newtheorem{theorem}{Theorem}

\theoremstyle{remark}

\newtheorem{remark}{Remark}
\theoremstyle{definition}

 \numberwithin{equation}{section}

\numberwithin{equation}{section}

\newcounter{comment}
\begin{document}

\title{On the connection problem for nonlinear differential equation}

\author{Zhao-Yun Zeng$^{}$\footnote{Corresponding author. Email:jenseng5@163.com.
}\ \ and ~Lin Hu}
 \date{\small School of Mathematics and Physics, Jinggangshan University, Ji'an 343009, PR China
}
\maketitle

\begin{abstract}
We consider the connection problem of the second nonlinear differential equation
\begin{equation}\label{nonlinear equ}
 \Phi''(x)=(\Phi'^2(x)-1)\cot\Phi(x)+ \frac{1}{x}(1-\Phi'(x))
\end{equation}
subject to the boundary condition $\Phi(x)=x-ax^2+O(x^3)$ ($a\geq0$) as $x\to0$. In view of that equation \eqref{nonlinear equ} is equivalent to the fifth Painlev\'{e} (PV) equation after a M\"{o}bius transformation, we are able to study the connection problem of equation \eqref{nonlinear equ} by investigating the corresponding connection problem of PV. Our research technique is based on the method of {\it uniform asymptotics} presented by Bassom el at. The monotonically solution on real axis of equation \eqref{nonlinear equ} is obtained, the explicit relation (connection formula) between the constants in the solution and the real number $a$  is also obtained. This connection formulas have been established earlier by Suleimanov via the isomonodromy deformation theory and the WKB method, and recently are applied for studying level spacing functions.
\end{abstract}


\vspace{5mm}


\section{Introduction and main results }
In the present paper we show how the technique of {\it uniform asymptotics} introduced by Bassom, Clarkson, Law and
McLeod in \cite{APC} can be applied to the equation
\begin{equation}\label{Phi equation}
  \Phi''(x)=(\Phi'^2(x)-1)\cot\Phi(x)+ \frac{1}{x}(1-\Phi'(x)),
\end{equation}
which solutions related to the computation of one particle density matrix of impenetrable bosons at zero temperature\cite{Creamer,Vaida:Tracy}.

We are focus on the problem of calculating an asymptotics behavior as $x\rightarrow\infty$ of one-parameter class of regular solutions
to equation \eqref{Phi equation} defined with the boundary condition
\begin{equation}\label{Expanding of Phi at 0}
  \Phi(x)=x-ax^2+O(x^3),\quad {\rm as}~x\to0,
\end{equation}
and on the relevant connection formulae between the different asymptotic parameters which appeared in the above mentioned critical expansions.

Introducing the change of variable
\begin{equation}\label{Trans:Phi:Pv}
  y(s)=\exp(-2i\Phi(x)),\quad s=\frac{x}{2},
\end{equation}
in Eq.\eqref{Phi equation}, we get for $y(s)$ the special fifth Painlev\'e (PV)
\begin{equation}\label{PV:RM}
   \frac{d^2y}{ds^2}=\Big(\frac{1}{2y}+\frac{1}{y-1}\Big)\Big(\frac{dy}{ds}\Big)^2-
\frac{1}{s}\frac{dy}{ds}-4i\frac{y}{s}+8\frac{y(y+1)}{y-1},
\end{equation}
which appears in the studying the level spacing functions related to the Fredholm determinant of the sine kernel $\frac{\sin\pi(x-y)}{\pi(x-y)}$
 on the finite interval $(-s,s)$ \cite{Shukla:1995,Novok2001}.
Let $s=it$ and $p(t)=\frac{\sqrt{y}+1}{\sqrt{y}-1}$, then equation \eqref{PV:RM} is equivalent to a special third Painlev\'{e} (PIII) equation
\begin{equation}\label{PIII:surface}
  \frac{d^2p}{dt^2}=\frac{1}{p}\Big(\frac{dp}{dt}\Big)^2-\frac{1}{t}\frac{dp}{dt}-\frac{1}{t}(p^2-1)+p^3-\frac{1}{p},
\end{equation}
which closely related to the studying Bonnet surfaces\cite{Bobe,Bobenko:book}, and the mean curvature and the metric in terms of $p(t)$ (see (3.115) in \cite{Bobenko:book}).
If we set $w(t)=-p(t)$, then $w(t)$ satisfies the another special PIII
\begin{equation}\label{PIII}
  \frac{d^2w}{dt^2}=\frac{1}{w}\Big(\frac{dw}{dt}\Big)^2-\frac{1}{t}\frac{dw}{dt}+\frac{1}{t}(w^2-1)+w^3-\frac{1}{w}.
\end{equation}
We mention that this special PIII can be expressed algebraically in terms of a fifth Painlev\'{e} transcendent and its first derivative.
 Consider the following pair of equations 
\begin{equation}\label{transf:PV:PIII}
 h(\tau)=\frac{w'(t)-w^2(t)-1}{w'(t)-w^2(t)+1},\qquad w(t)=\frac{2\tau h(\tau)}{\tau h'(\tau)-h(\tau)+1},
\end{equation}
where $\tau=\frac{t^2}{2}$. Eliminating $h$ from \eqref{transf:PV:PIII}, we get Eq.\eqref{PIII} for $w(t)$;
and eliminating $w$ from \eqref{transf:PV:PIII}, we get for $h(\tau)$ the special PV equation
\begin{equation}\label{PV}
 \frac{d^2h}{d\tau^2}=(\frac{1}{2h}+\frac{1}{h-1})\Big(\frac{dh}{d\tau}\Big)^2-\frac{1}{\tau}\frac{dh}{d\tau}
 -\frac{1}{8}\frac{(h-1)^2}{\tau^2h}-\frac{h}{\tau},
\end{equation}
which admits the Lax representation \cite{Sulei}.

With the help of the preceding derivation we can now see that the Eq.\eqref{Phi equation} is equivalent to Eq.\eqref{PV} after the transform
\begin{equation}\label{relation:h:Phi}
  h(\tau)=1+\frac{2\sin^2\frac{\Phi(x)}{2}}{\Phi'(x)-1},\quad \tau=-\frac{x^2}{8}.
\end{equation}
Hence, we are able to study the connection problem of Eq.\eqref{Phi equation} by using the isomonodromic deformation technique \cite{Its book} to
consider the correspond connection problem for PV \eqref{PV}. Basing on the special Lax pair of \eqref{PV},
the author of \cite{Sulei} has studied analytical this solution to Eq.\eqref{Phi equation} with the initial condition \eqref{Expanding of Phi at 0}
for all real-valued $a$,
and obtained the asymptotic expansion of $\Phi$ as $x\to\infty$ and explicit connection formulas by virtue of the isomonodromic deformation
technique. Let us summarize the main result of \cite{Sulei} in the theorem as follows:
\begin{theorem}\label{theo:solu:conect}
There exists a unique solution of \eqref{Phi equation} which satisfies \eqref{Expanding of Phi at 0} for any given real number $a$.
\begin{itemize}
  \item [(A)] If $a>\frac{1}{\pi}$, this solution  exists for positive real $x$, and
  \begin{equation}\label{solution:3}
    \Phi(x,a)=-x+\beta\ln x +\gamma +o(1),~{\rm as}~ x\to\infty,
  \end{equation}
  where $\beta$ and $\gamma$ are real constants. Furthermore, the relationship between the parameters $\beta,~\gamma$ in \eqref{solution:3} and the parameters
  $a$ in \eqref{Expanding of Phi at 0} are provided by the connection formulas
  \begin{align}
    \beta &=-\frac{1}{\pi}\ln(a\pi-1), \label{connect:formula:2}\\
         \gamma &=\frac{\pi}{2}+2\arg\Gamma(\frac{i\beta}{2}-\frac{1}{2})+\beta\ln2+k\pi,~~k\in\mathbb{Z}. \label{connect:formula:2:2}
  \end{align}
  \item [(B)] If $a<\frac{1}{\pi}$, this solution  exists for all real $x$, and increases monotonically as $x\to\infty$,
  \begin{equation}\label{solution:1}
    \Phi(x,a)=x+\beta\ln x +\gamma +o(1),
  \end{equation}
  where $\beta$ and $\gamma$ are real constants. Furthermore, the relationship between the parameters $\beta,~\gamma$ in \eqref{solution:1} and the parameters
  $a$ in \eqref{Expanding of Phi at 0} are provided by the connection formulas
  \begin{align}
       \beta &=\frac{1}{\pi}\ln(1-a\pi), \label{connect:formula:1}\\
         \gamma &=-2\arg\Gamma(\frac{i\beta}{2})+\beta\ln2-\pi{\rm sign}\beta, \label{connect:formula:1:2}
  \end{align}
  where $\beta\neq0$ and $\gamma(0)=0$.
  \item [(C)] If $a=\frac{1}{\pi}$, this solution  exists for all real $x$ as $x$ increases monotonically  to a finite limit, and, as $x\to\infty$,
  \begin{equation}\label{solution:2}
    \Phi(x,\frac{1}{\pi})=\frac{\pi}{2} +o(1).
  \end{equation}
\end{itemize}
\end{theorem}
On the earlier work of the authors of \cite{Creamer}, they have studied numerically this solution with the given asymptotic behavior
 at the origin \eqref{Expanding of Phi at 0} for the case of $a>\frac{1}{\pi}$ and proposed \eqref{solution:3} and \eqref{connect:formula:2}, however,
  they not obtain the explicit expression \eqref{connect:formula:2:2} for $\gamma$. Recently, the connection formulas in
  Theorem \ref{theo:solu:conect} are
applied for calculations of the Fredholm determinant of the sine kernel $\sin\pi(x-y)/\pi(x-y)$ on the finite interval $(t,-t)$; see \cite{Novok2001}.

In this paper, we provided a simpler and more rigorous proof of the Theorem \ref{theo:solu:conect}, by using the {\it uniform asymptotics} method proposed in \cite{APC}. 
For our purposes, we first briefly outline some important properties of the theory of monodromy preserving deformations for the PV transcendents. The reader is referred to \cite{Fokas book, MJ} for more details.

One of the Lax pairs for the fifth Painlev\'e equation \eqref{PV} is the system of linear ordinary equations \cite{Novok2001}
\begin{equation}\label{PV:lax pair 1}
\frac{\partial\Psi}{\partial\lambda}=\left\{
-i\tau\sigma_3+
\frac{1}{\lambda}\left(\begin{array}{cc}
\frac{1}{4}&u\\
v&-\frac{1}{4}
\end{array}\right)
+\frac{1}{\lambda^2}\left(
                    \begin{array}{cc}
                      z & q \\
                      q & -z \\
                    \end{array}
                  \right)
\right\}\Psi
\end{equation}
and
\begin{equation}\label{PV:lax pair 2}
  \frac{\partial\Psi}{\partial \tau}=\left\{
 -i\lambda\sigma_3 +\left(
                      \begin{array}{cc}
                        g              & \frac{u}{\tau} \\
                        \frac{v}{\tau} & -g             \\
                      \end{array}
                    \right)
\right\}\Psi,
\end{equation}
where $\tau=-\frac{x^2}{8}$ and
\begin{equation}
z=-\frac{i}{8}\frac{h+1}{h-1},\, q=-\frac{1}{4}\frac{\sqrt{h}}{h-1},\,
u+v=-\frac{i}{2\sqrt{h}},\, u-v=\frac{i\tau h_\tau}{(1-h)\sqrt{h}},\,g=\frac{1}{8\tau}(1+\frac{1}{h}).\label{PV:u:v}
 \end{equation}
The compatibility condition $\Psi_{\lambda x}=\Psi_{x\lambda}$ implies that $h(\tau)$ satifies the PV equation \eqref{PV}.

The equation \eqref{PV:lax pair 1} has two irregular points $\lambda=0$ and $\lambda=\infty$. There exists a canonical solution $\Psi^{(\infty)}$
defined in a neighborhood of the irregular singular point $\lambda=\infty$ with the following asymptotics behavior
\begin{equation} \label{Psi boundary at infty1}
  \Psi^{(\infty)}(\lambda)=E^{(\infty)}(\tau)(I+O(\frac{1}{\lambda}))\lambda^{\frac{1}{4}\sigma_3}\exp(-i\tau\lambda\sigma_3),
  ~\lambda\rightarrow\infty,~\arg\lambda=0,
 \end{equation}
 where
 \begin{equation}\label{Express E:infty}
  E^{(\infty)}(\tau)=\tau^{\frac{1}{8}\sigma_3}\exp(\sigma_3 J(\tau)):=d^{\sigma_3}
 \end{equation}
 with $J(\tau)=\frac{1}{8}\int_{-c}^{\tau}\frac{dt}{th(t)}$, here $c$ is a postive constant.

 From \eqref{PV:lax pair 1} and \eqref{Psi boundary at infty1} it follows that
 \begin{equation}\label{Psi boundary at infty2}
   \Psi^{(\infty)}(\lambda)=E^{(\infty)}(\tau)\left(
                                                \begin{array}{cc}
                                                  1+O(\frac{1}{\lambda}) & \frac{a_1}{\lambda} \\
                                                  \frac{a_2}{\lambda} & 1+O(\frac{1}{\lambda}) \\
                                                \end{array}
                                              \right)\lambda^{\frac{1}{4}\sigma_3}\exp(-i\tau\lambda\sigma_3)
 \end{equation}
  with $a_1=\frac{u}{2i\tau d}$ and $a_2=-\frac{v}{2i\tau}d$.

 On the other hand, \eqref{PV:lax pair 1} has another canonical solution $\Psi^{(0)}$ in a neighbourhood of the irregular singular point $\lambda=0$
 \begin{equation}\label{Psi boundary at zero}
  \Psi^{(0)}(\lambda)=H(\tau)E^{(0)}(\tau)(I+O(\lambda))\lambda^{\frac{1}{4}\sigma_3}\exp(\frac{i}{\lambda}\sigma_3),
  ~\lambda\rightarrow0,~\arg\lambda=0,
 \end{equation}
 where the coefficients $H$ and $E^{(0)}$ have the form
 \begin{equation}\label{matrice H:E}
   H(\tau)=\frac{1}{\sqrt{h(\tau)-1}}(i\sigma_3\sqrt{h}+\sigma_1),
   \quad E^{(0)}(\tau)=\tau^{\frac{1}{8}\sigma_3}\exp(-\sigma_3 J(\tau))=\tilde{d}^{\sigma_3}.
 \end{equation}
 Since $\Psi^{(\infty)}$ and $\Psi^{(0)}$ are both fundamental solutions, the connection matrix $Q$ can be defined by
 \begin{equation}\label{Connect matrice}
   \Psi^{(\infty)}(\lambda)=\Psi^{(0)}(\lambda)Q,
 \end{equation}
  Differentiating both sides
of equation \eqref{Connect matrice} with respect to $x$, and making use of the fact that both $\Psi^{(\infty)}$ and $\Psi^{(0)}$ satisfy \eqref{PV:lax pair 1}, it is easily found the isomonodromic condition $\frac{dQ}{dx}=0$; i.e., $Q$ is a constant matrix.

In the framework of the isomonodromic deformation method, one needs to calculate the monodromy data $Q$ both in terms of the initial condition \eqref{Expanding of Phi at 0} and asymptotics \eqref{solution:1}. Equating then the leading terms of nontrivial monodromy data, one gets
connection formulas for the parameters $\beta$, $\gamma$ and $a$. In the limit $x\to0$ the first term of equation in $\lambda$ \eqref{PV:lax pair 1} vanishes, so the $\Psi$-functions can be expressed via the Whittaker
 functions \cite{Sulei}, then the monodromy data as $x\to0$ is calculated explicitly by use of multiplication formulas of the Whittaker functions 
 and obtain that 
\begin{align}
  (Q)_{21}&=i2^{-3/4}\sqrt{a\pi}.\label{Q:0}
\end{align}
To estimate the connection matrix in the limit $\tau\to+\infty$. One finds the WKB solution of the $\Psi$-function, which demands a standard procedure of matching near the turning points, involving parabolic cylinder functions(see,\cite[p.251]{Sulei}). Eventually, one can obtains the connection matrix for large $x$ as follows by using of the asymptotics behavior of parabolic cylinder functions.

In this paper, we shall provide a hopefully simpler and more rigorous derivation of the asymptotic behavior and the connection formulas in Theorem \ref{theo:solu:conect} by using
the {\it uniform asymptotics} method presented in \cite{APC}. Along the same lines we may find the work of Olver \cite{FO} and Dunster  \cite{Dunster} for coalescing turning points.
Initially in \cite{APC}, the second Painlev\'e(PII) equation has been taken as an example to illustrate the method. While the difficulty in extending the techniques for PII to other
transcendents is also acknowledged by the authors of \cite[p.244]{APC}.
Yet the  method has been applied to the connection problems
by Wong and Zhang \cite{WZ1,WZ2}, Zeng and Zhao\cite{Zeng}, Long, Zeng and Zhou \cite{LongZ}.
Recently, Long el at.\cite{LongL} present a detail asymptotics analysis of the real solutions of the first Painlev\'e(PI) equation by virtue of the {\it uniform asymptotics} method.

 The rest of the paper is organized as follows. The proof of Theorem \ref{theo:solu:conect} is provided in the Section \ref{Proof theorem}. In Section \ref{Uniform infty}, for the case of $a>\frac{1}{\pi}$ and $a<\frac{1}{\pi}$, we derive  uniform approximations to the solutions of the second-order differential equation obtained from the Lax pair  \eqref{PV:lax pair 1}  as $x\rightarrow+\infty$ by virtue of the parabolic cylinder functions on the Stokes curves, respectively. The entry $(2,1)$ of the connection matrix $Q$ for large $x$ is also computed in the section. Some technical details are put in Appendices A and B to clarify the derivation.
\section{Proof of the theorem}\label{Proof theorem}
To proof of the Theorem \ref{theo:solu:conect}, we need two lemmas as follows.
\begin{lemma}\label{Q:a>pi}
For $a>\frac{1}{\pi}$, the asymptotics behavior of the entry $(2,1)$ of the connection matrix $Q$ is
\begin{equation}\label{Q1:-infty}
(Q)_{21}=\frac{2^{-\frac{1}{4}}\sqrt{\pi}{\rm e}^{-\frac{\pi\beta }{4}}}{\Gamma(\frac{1}{2}-i\frac{\beta}{2})}\exp\Big(iS+\frac{i}{2}x-\frac{i\beta}{2}\ln x-\frac{i\beta}{2}\ln2
  +\frac{3\pi i}{4}\Big).
\end{equation}
\end{lemma}
\begin{lemma}\label{Q:a<pi}
For $a<\frac{1}{\pi}$, the asymptotics behavior of the entry $(2,1)$ of the connection matrix $Q$ is
\begin{equation}\label{Q2:-infty}
  (Q)_{21}=\frac{i\sqrt{\beta}2^{-\frac{3}{4}}\sqrt{\pi}{\rm e}^{\frac{\pi\beta }{4}}}{\Gamma(\frac{i\beta}{2}+1)}
  \exp\Big(-iS+\frac{i}{2}x+\frac{i\beta}{2}\ln x+\frac{i\beta}{2}\ln2\Big).
\end{equation}
\end{lemma}
The rigorous proofs of those results will be given in the next section. With the help of the preceding two lemmas we can now prove Theorem \ref{theo:solu:conect}.

\noindent {\bf{Proof of Theorem \ref{theo:solu:conect}}}: We first give the proof when $a>\frac{1}{\pi}$. Since the connection matrix $Q$ must be independent of $x$, it follows that the right hand sides of \eqref{Q:0} and \eqref{Q1:-infty} are equality. Separating real and imaginary parts and in view of the standard formulae $|\Gamma(\frac{1}{2}+iy)|^2=\frac{\pi}{\cosh y\pi}$(see,\cite{OL}) gives  the asymptotic behaviors \eqref{solution:3} and connection formulas \eqref{connect:formula:2},\eqref{connect:formula:2:2}, which completes the proof of statement (A) in Theorem \ref{theo:solu:conect}.

When $a<\frac{1}{\pi}$, the asymptotic behaviors \eqref{solution:1} and connection formulas \eqref{connect:formula:1},\eqref{connect:formula:1:2} are obtained straightforward by equating the expressions \eqref{Q:0} and \eqref{Q2:-infty}, here has been made of the formulas
$\Gamma(z+1)=z\Gamma(z)$ and $|\Gamma(iy)|^2=\frac{2\pi}{y({\rm e}^{y\pi}-{\rm e}^{-y\pi})}$ and the fact that $\Phi(x,0)=x$ is the solution of the initial problem \eqref{Phi equation} and \eqref{Expanding of Phi at 0}. Hence the statement (B) in Theorem \ref{theo:solu:conect} is proved.

Using the fact that $\Phi(x,a)$ is a continuous function of $a$ \cite[Lemma1,p253]{Sulei}, and taking \eqref{solution:1} into consideration, according to the definition of $L$ in \eqref{Defin L}, we obtain
\begin{equation*}
\lim_{x\to\infty}\Phi(x,\frac{1}{\pi})=\frac{\pi}{2},
\end{equation*}
which gives the proof of statement (C) in Theorem \ref{theo:solu:conect}.

The proof of Theorem \ref{theo:solu:conect} is now complete.\qed
\section{Uniform asymptotics and proofs of the lemmas }\label{Uniform infty}
Making the scaling
\begin{equation}\label{scale transform}
  \xi=x,\quad\eta=x\lambda,
\end{equation}
so that \eqref{PV:lax pair 1} becomes
\begin{equation}\label{PV:lax pair new}
  \frac{\partial\Psi}{\partial\eta}=\left(
                                      \begin{array}{cc}
                                        \frac{i}{8}\xi+\frac{1}{4\eta}+\frac{\xi}{\eta^2}z & \frac{u}{\eta}+\frac{\xi}{\eta^2}q \\
                                        \frac{v}{\eta}+\frac{\xi}{\eta^2}q & -(\frac{i}{8}\xi+\frac{1}{4\eta}+\frac{\xi}{\eta^2}z) \\
                                      \end{array}
                                    \right)\Psi
\end{equation}
Let $(\Psi_1,\Psi_2)^T$ be a independent solution of \eqref{PV:lax pair new}, and set
\begin{equation}\label{define phi}
  \phi=\Big(\frac{v}{\eta}+\frac{\xi}{\eta^2}q\Big)^{-\frac{1}{2}}\Psi_2,
\end{equation}
we get from \eqref{PV:lax pair new} the second-order linear differential equation for $\phi(\eta)$
\begin{align}\label{Second order equation}
  \frac{d^2\phi}{d\eta^2}
  &=\Big\{\xi^2\Big(\frac{i}{8}+\frac{z}{\eta^2}\Big)^2+\frac{\xi}{2\eta}(\frac{i}{8}+\frac{z}{\eta^2})+\frac{1}{16\eta^2}+\xi^2\frac{q^2}{\eta^4}
  +\frac{1}{\eta^2}[uv+\frac{\xi}{\eta}q(u+v)]\nonumber\\
  &\quad+\frac{1}{4\eta^2}+\xi\frac{2z}{\eta^3}-\xi\frac{1}{\eta}(\frac{i}{8}+\frac{z}{\eta^2})l_1-\frac{1}{4\eta^2}\l_1
  +\frac{3}{4\eta^2}l_1^2-\frac{1}{\eta^2}l_2
  \Big\}\phi,
\end{align}
where
\begin{equation}\label{define l1}
 l_1=\frac{v+\frac{2\xi}{\eta}q}{v+\frac{\xi}{\eta}q},\quad
l_2=\frac{v+\frac{3\xi}{\eta}q}{v+\frac{\xi}{\eta}q}.
\end{equation}
From \eqref{PV:u:v} it is easy to verify that
\begin{equation}\label{z2:q2}
  z^2+q^2=-\frac{1}{64},\quad q(u+v)=\frac{i}{8(h-1)}.
\end{equation}
Substituting \eqref{z2:q2} into \eqref{Second order equation} yields
\begin{align}\label{Second order equation2}
  \frac{d^2\phi}{d\eta^2}&=\Big\{-\frac{\xi^2}{64}(1-\frac{1}{\eta^2})^2+\frac{\xi^2}{64\eta^2}(16iz-2)+\frac{\xi}{2\eta}(\frac{i}{8}+\frac{z}{\eta^2})
  -\frac{\xi}{\eta}(\frac{i}{8}+\frac{z}{\eta^2})l_1\nonumber\\
  &\quad+\frac{uv}{\eta^2}+\xi\frac{2z}{\eta^3}+\frac{i\xi}{8\eta^3(h-1)}+\frac{1}{\eta^2}(\frac{5}{16}-\frac{1}{4}l_1+\frac{3}{4}l_1^2-\l_2)
  \Big\}\phi.
\end{align}
\subsection{Proof of the Lemma \ref{Q:a>pi}}
To proof the Lemma \ref{Q:a>pi}, we need several lemmas. First, we need to construct the uniform asymptotics solution of equation \eqref{Second order equation2} as $\xi\rightarrow\infty$ for $a>\frac{1}{\pi}$. When $a>\frac{1}{\pi}$, according to \cite[(A.71)]{Sulei}, the solution of boundary value problem \eqref{Phi equation}-\eqref{Expanding of Phi at 0} has the following asymptotics expansion
\begin{equation}\label{asy derivative of Phi}
  \Phi'(\xi,a)=-1+\frac{\varphi(S)}{\xi}+O(\xi^{-2}),~{\rm as}~\xi\to\infty,
\end{equation}
where $\varphi(S)=\sin 4S+2k^2\sin^22S$ with $S=\frac{1}{2}\Phi(\xi,a)$.
It follows from the expression of $h$ in \eqref{relation:h:Phi} that
\begin{equation}\label{asy h}
  h(\xi)
  =\cos^2S\Big(1-\frac{\varphi(S)\tan^2S}{2\xi}+O(\xi^{-2})\Big),~as~ \xi\to\infty.
\end{equation}
From \eqref{asy h} and \eqref{PV:u:v}, we obtain the following asymptotic behaviors as $\xi\to\infty$
\begin{align}
  & z=\frac{i}{8}\Big(1+2\cot^2S-\frac{\varphi(S)\csc^2S}{\xi}+O(\xi^{-2})\Big),\label{asym z}\\
  &uv=\frac{\xi^2}{16\sin^2S}\Big(1-\frac{2\tan S+\varphi(S)\big(1+\frac1 2 \tan^2S\big)}{\xi}+O(\xi^{-2})\Big), \label{asym uv}\\
  &\frac{q}{v}
  =\frac{i\cot S}{\xi}(1+O(\xi^{-1})),\label{asym q:v}\\
  &qv=-\frac{i\xi\cot S}{16\sin^2S}\Big(1+O(\xi^{-1})\Big).\label{asym qv}
\end{align}
Then, for large $\xi$, substituting \eqref{asym z}, \eqref{asym uv} and \eqref{asym q:v} into \eqref{Second order equation2}, a tedious but straightforward calculation gives
\begin{equation}\label{Second order equation asy}
  \frac{d^2\phi}{d\eta^2}
  =-\xi^2F(\xi,\eta)\phi,
\end{equation}
where
\begin{equation}\label{expr:F}
  F(\xi,\eta)=\frac{1}{64}\Big(1-\frac{1}{\eta^2}\Big)^2+\frac{F_1(\xi,\eta)}{\xi}+F_2(\eta)O(\frac{1}{\xi^2}),
\end{equation}
here
\begin{equation}\label{express:F:1}
  F_1(\xi,\eta)=\frac{k^2}{4\eta^2}+\frac{i}{8\eta}\Big(1-\frac{1}{\eta^2}\Big)\Big(\frac{1}{2}+\frac{1}{b\eta-1}\Big)
\end{equation}
with $b=i\tan S$, and
\begin{equation}\label{express:F:2}
 F_2(\eta)=\frac{1}{\eta^2}+\frac{1}{\eta^3}.
\end{equation}

For large $\xi$, it follows from equation \eqref{Second order equation asy} that there are two coalescing turning
points near $\eta=1$, and two close to $\eta=-1$. In the present paper, we are only concerned with the two turning points, say $\eta_1$ and $\eta_2$,
near $\eta=1$. When $\eta_j$ approach to 1, it follows from \eqref{express:F:1}  that
\begin{equation}\label{asy F:1}
  F_1(\xi,1)\thicksim\frac{k^2}{4}.
\end{equation}
By using  \eqref{Second order equation asy} and \eqref{asy F:1}, we get the asymptotic formulas for the two turning points
\begin{equation}\label{Two turning points}
  \eta_j^{-1}=1\pm 2\xi^{-1/2}\sqrt{k^2}(1+o(1)),~j=1,2,
\end{equation}
 which  coalescing to $1$ when $\xi\rightarrow\infty$, and the Stokes' curves defined by
\begin{equation}\label{Stokes curves}
  \Im(\xi(\eta+\frac{1}{\eta}))=0.
\end{equation}
Assuming that $\xi\in\mathbb{R}^+$, then, it follows from \eqref{Stokes curves}  that the Stokes lines of the solution $\phi$ to
\eqref{Second order equation asy} are the positive and the negative real lines in the $\eta$ plane.

According to the  philosophy of uniform asymptotics in \cite{APC}, we define a number $\alpha$ by
\begin{equation}\label{Define alpha}
  \frac{1}{2}\pi i\alpha^2=\int_{-\alpha}^\alpha(\tau^2-\alpha^2)^{1/2}d\tau=\int_{\eta_1}^{\eta_2}F^{1/2}(\xi,s)d s,
\end{equation}
and a new variable $\zeta$ by
\begin{equation}\label{Define zeta}
  \int_{\alpha}^{\zeta}(\tau^2-\alpha^2)^{1/2}d\tau=\int_{\eta_2}^\eta F^{1/2}(\xi,s)ds.
\end{equation}
Here and in \eqref{Define alpha},
the cut for the integrand on the left-hand side  is the line segment joining $-\alpha$ and  $\alpha$. The path of integration is taken   along the upper edge of the cut.
With $\alpha$ and $\zeta$ so chosen, then the following lemma is a result from \cite[Theorem 1]{APC}.
\begin{lemma}\label{Uniform asymptotic theorem}~Given any solution $\phi(\eta,\xi)$ of \eqref{Second order equation asy},
there exist constants $c_1,~c_2$ such that, uniformly for $\eta$ on the Stokes cures defined by $(\ref{Stokes curves})$, as $\xi\rightarrow+\infty$,
 \begin{equation}\label{Uniform approximate of phi}
   \phi(\eta,\xi)=\Big(\frac{\zeta^2-\alpha^2}{F(\xi,\eta)}\Big)^{\frac{1}{4}}
\Big\{[c_1+o(1)]D_{\nu}({\rm e}^{\pi i/4}\sqrt{2\xi}\zeta)
+[c_2+o(1)]D_{-\nu-1}({\rm e}^{-\pi i/4}\sqrt{2\xi}\zeta)\Big\},
 \end{equation}
 where $D_{\nu}(z)$ and $D_{-\nu-1}(z)$ are solutions of the parabolic cylinder equation and $\nu$ defined by
\begin{equation}\label{Define nu}
  \nu=-\frac{1}{2}+\frac{1}{2}i\xi\alpha^2.
\end{equation}
\end{lemma}
The next thing to do in calculating the connection martix $Q$ as $\xi\rightarrow+\infty$ is to clarify the relation  between $\zeta$ and $\eta$ in \eqref{Define zeta}.
\begin{lemma}\label{lem relation zeta and eta infty}
For large $\xi$ and $\eta$,
\begin{equation}\label{relation zeta and eta infty}
  \frac{1}{2}\zeta^2=\frac{\alpha^2}{2}\ln\zeta+\frac{1}{8}(\eta+\frac{1}{\eta})-\frac{1}{4}+\frac{i}{4\xi}\ln\eta - \frac{i}{2\xi}\ln(1-b^{-1})+o(\xi^{-1}),
\end{equation}
where $b=i\tan S$, and
\begin{equation}\label{asy alpha}
  \alpha^2=-\frac{k^2}{\xi}+o(\frac{1}{\xi})\quad as~\xi\to\infty.
\end{equation}
\end{lemma}
\begin{remark}
Coupling \eqref{Define nu} and \eqref{asy alpha} determines the approximate value
\begin{equation}\label{asym:nu}
  \nu=-\frac{ik^2}{2}-\frac 1 2+o(1)\quad {\rm as}~\xi\to+\infty
\end{equation}
for the order of the parabolic cylinder function $D_{\nu}({\rm e}^{\pi i/4}\sqrt{2\xi}\zeta)$ in \eqref{Uniform approximate of phi}.
\end{remark}
\begin{lemma}\label{lem relation zeta and eta zero}
When $\eta\to 0$, for large $\xi$, such $\xi\eta=o(1)$, then holds
\begin{equation}\label{relation zeta and eta zero}
  \frac{1}{2}\zeta^2=\frac{\alpha^2}{2}\ln\zeta+\frac{1}{8}(\eta+\frac{1}{\eta})-\frac{1}{4}-\frac{i}{4\xi}\ln\eta - \frac{i}{2\xi}\ln(1-b)-\frac{1}{2}\pi i\alpha^2+o(\xi^{-1}).
\end{equation}
\end{lemma}

The proofs of Lemmas \ref{lem relation zeta and eta infty} and \ref{lem relation zeta and eta zero} will be given in Appendix A and B, respectively. We are now turning to the proof of Lemma \ref{Q:a>pi}.

\noindent {\bf{Proof of Lemma \ref{Q:a>pi}}}.
We will concentrate on evaluating the connection matrice $Q$, for which we need the uniform
asymptotic behaviors of $\phi$. From the definition of $Q$ \eqref{Connect matrice} it follows that
\begin{align}
\Psi_{21}^{(\infty)}=(Q)_{11}\Psi_{21}^{(0)}+(Q)_{21}\Psi_{22}^{(0)}.\label{Express Phi:21:infty:0}
\end{align}
Hence, the first task is to find the expressions of $\Psi_{21}^{(\infty)},~\Psi_{21}^{(0)}$ and $\Psi_{22}^{(0)}$, respectively.

For large $\xi$, it follows from lemma \ref{Uniform asymptotic theorem} that
 two linearly independent asymptotic solutions of equation \eqref{Second order equation asy} are $\tilde{\phi}_{\nu}$ and $\tilde{\phi}_{-\nu-1}$
which are uniform with respect to $\eta$ on the Stokes' curves. Here
\begin{equation}\label{Denote of tilde phi1}
  \tilde{\phi}_{\nu}=\Big(\frac{\zeta^2-\alpha^2}{F(\xi,\eta)}\Big)^{\frac{1}{4}}D_{\nu}({\rm e}^{\pi i/4}\sqrt{2\xi}\zeta)
\end{equation}
and
\begin{equation}\label{Denote of tilde phi2}
  \tilde{\phi}_{-\nu-1}=\Big(\frac{\zeta^2-\alpha^2}{F(\xi,\eta)}\Big)^{\frac{1}{4}}D_{-\nu-1}({\rm e}^{-\pi i/4}\sqrt{2\xi}\zeta).
\end{equation}

By virtue of \eqref{define phi}, we have
\begin{equation}\label{asy:Phi:21:infty}
\Psi_{21}^{(\infty)}= \Big(\frac{v}{\eta}+\frac{\xi q}{\eta^2})^{1/2}(f_1\tilde{\phi}_{\nu}+f_2\tilde{\phi}_{-\nu-1}),
\end{equation}
where $f_j(j=1,2)$ are undetermined constants which can be determined by \eqref{Psi boundary at infty2}.

Similarly, we obtain
\begin{align}
\Psi_{21}^{(0)}= \Big(\frac{v}{\eta}+\frac{\xi q}{\eta^2})^{1/2}(\delta_1\tilde{\phi}_{\nu}+\delta_2\tilde{\phi}_{-\nu-1}),\label{asy:Phi:21:orgin}\\
\Psi_{22}^{(0)}= \Big(\frac{v}{\eta}+\frac{\xi q}{\eta^2})^{1/2}(\delta_3\tilde{\phi}_{\nu}+\delta_4\tilde{\phi}_{-\nu-1}),\label{asy:Phi:22:orgin}
\end{align}
where $\delta_j(j=1,2,3,4)$ are undetermined constants which can be determined by \eqref{Psi boundary at zero}.

Substituting \eqref{asy:Phi:21:infty}, \eqref{asy:Phi:21:orgin} and \eqref{asy:Phi:22:orgin} into \eqref{Express Phi:21:infty:0}, then comparing with the
coefficients of $\tilde{\phi}_{\nu}$ and $\tilde{\phi}_{-\nu-1}$, respectively, we get
\begin{equation}\nonumber
  \left\{\begin{array}{ll}
    \delta_1(Q)_{11}+\delta_3(Q)_{21}=f_1 \\
    \delta_2(Q)_{11}+\delta_4(Q)_{21}=f_2
  \end{array}\right.
\end{equation}
which gives us that
\begin{equation}\label{express:Q:21}
  (Q)_{21}=\frac{\delta_1f_2-\delta_2f_1}{\delta_1\delta_4-\delta_2\delta_3}.
\end{equation}

To calculate $f_1,~f_2$ and $\delta_j~(j=1,\ldots,4)$ we proceed as follows.
We shall be interested in finding the asymptotic behavior of $\tilde{\phi}_{\nu}$ and $\tilde{\phi}_{-\nu-1}$ in \eqref{Denote of tilde phi1}
\eqref{Denote of tilde phi2} for $\eta$ on the Stokes line $\arg\eta=0$ as $\eta\to\infty$ and $\eta\to0$, respectively. Then, substituting the
obtained results into \eqref{asy:Phi:21:infty}, \eqref{asy:Phi:21:orgin} and \eqref{asy:Phi:22:orgin}, we will obtain the asymptotic behavior of
$\Psi_{21}^{(\infty)}$, $\Psi_{21}^{(0)}$ and $\Psi_{22}^{(0)}$, respectively, which contain the constants $f_1,~f_2$ and $\delta_j~(j=1,\ldots,4)$.
Combining with the boundary conditions \eqref{Psi boundary at infty2} and \eqref{Psi boundary at zero} for $\Psi^{(\infty)}$ and $\Psi^{(0)}$, one can
determine the constants $f_1,~f_2$ and $\delta_j~(j=1,\ldots,4)$.

From \cite{OL}, we have the asymptotic behavior of $D_{\nu}(z)$ for $|z|\rightarrow\infty$ as follows:
\begin{equation}\label{Asymptotic behavior of D z}
  D_{\nu}(z)\sim\left\{\begin{array}{ll}
z^{\nu}{\rm e}^{-\frac{1}{4}z^2},&\arg z\in(-\frac{3}{4}\pi,\frac{3}{4}\pi),\\
z^{\nu}{\rm e}^{-\frac{1}{4}z^2}-\frac{\sqrt{2\pi}}{\Gamma(-\nu)}{\rm e}^{i\pi\nu}z^{-\nu-1}{\rm e}^{\frac{1}{4}z^2},&\arg z=\frac{3}{4}\pi,\\
{\rm e}^{-2\pi i\nu}z^{\nu}{\rm e}^{-\frac{1}{4}z^2}-\frac{\sqrt{2\pi}}{\Gamma(-\nu)}{\rm e}^{i\pi\nu}z^{-\nu-1}{\rm e}^{\frac{1}{4}z^2},
&\arg z=\frac{5}{4}\pi.
 \end{array}
 \right.
\end{equation}
For $\eta$ on the Stokes line $\arg\eta=0$ and $\eta\rightarrow\infty$,
From \eqref{relation zeta and eta infty} it immediately follows
that $\zeta^2\thicksim \frac{1}{4}\eta$, then we have $\arg\zeta\sim0$. Therefore,
$\arg({\rm e}^{\pi i/4}\sqrt{2\xi}\zeta)\thicksim\frac{\pi}{4}$ and $\arg({\rm e}^{-\pi i/4}\sqrt{2\xi}\zeta)\thicksim-\frac{\pi}{4}$ for $\xi>0$.
From \eqref{Second order equation asy} we have $F^{-1/4}\thicksim 2^{\frac{3}{2}}$ as $\eta\rightarrow\infty$ for large $\xi$. Since
$(\zeta^2-\alpha^2)^{1/4}\thicksim\zeta^{1/2}$ as $\eta\rightarrow\infty$, by using the appropriate asymptotic formulas of $D_{\nu}(z)$ in
(\ref{Asymptotic behavior of D z}), we obtain from \eqref{Denote of tilde phi1} and \eqref{relation zeta and eta infty} that
\begin{align}\label{asy tilde phi:1}
\tilde{\phi}_{\nu}
\thicksim A_0\eta^{\frac{1}{4}}{\rm e}^{-\frac{1}{8}i\xi\eta},~{\rm as}~\eta\rightarrow\infty,
\end{align}
here has been made of \eqref{Define nu}, where
\begin{equation}\label{coef:A0}
  A_0=2^{\frac{\nu+3}{2}}{\rm e}^{\frac{\pi i}{4}\nu} {\rm e}^{\frac{i}{4}\xi+\frac{\nu}{2}\ln\xi}\big(1-b^{-1}\big)^{-\frac{1}{2}}.
\end{equation}
Similarly, from \eqref{Denote of tilde phi2} and \eqref{relation zeta and eta infty} it follows that
\begin{align}\label{asy tilde phi:2}
\tilde{\phi}_{-\nu-1}
\thicksim B_0\eta^{-\frac{1}{4}}{\rm e}^{\frac{1}{8}i\xi\eta},~{\rm as}~\eta\rightarrow\infty,
\end{align}
where
\begin{equation}\label{coef:B0}
  B_0=2^{1-\frac{\nu}{2}}{\rm e}^{\frac{\pi i}{4}(\nu+1)} {\rm e}^{-\frac{i}{4}\xi-\frac{\nu+1}{2}\ln\xi}\big(1-b^{-1}\big)^{\frac{1}{2}}.
\end{equation}
Substituting \eqref{asy tilde phi:1} and \eqref{asy tilde phi:2} into \eqref{asy:Phi:21:infty}, yields
\begin{equation}\label{asy:behavior:Phi:21:infty}
  \Psi_{21}^{(\infty)}(\eta)\thicksim v^{\frac{1}{2}}\eta^{-\frac{1}{2}}\Big(f_1A_0\eta^{\frac{1}{4}}{\rm e}^{-\frac{1}{8}i\xi\eta}
  +f_2B_0\eta^{-\frac{1}{4}}{\rm e}^{\frac{1}{8}i\xi\eta}\Big)~{\rm as}~\eta\rightarrow\infty.
\end{equation}
Moreover, it has the asymptotic behavior prescribed in \eqref{Psi boundary at infty2} when $\eta\to\infty$. Thus, we have
\begin{equation}\label{matching:infty}
 v^{\frac{1}{2}}\eta^{-\frac{1}{2}}\Big(f_1A_0\eta^{\frac{1}{4}}{\rm e}^{-\frac{1}{8}i\xi\eta}
  +f_2B_0\eta^{-\frac{1}{4}}{\rm e}^{\frac{1}{8}i\xi\eta}\Big)\thicksim -\frac{d}{2i\tau}v\lambda^{-\frac{3}{4}}\exp(-i\tau\lambda)
\end{equation}
Comparing the coefficients of ${\rm e}^{\frac{1}{8}i\xi\eta}$ and ${\rm e}^{-\frac{1}{8}i\xi\eta}$ on both sides of the above asymptotic equation, we get
\begin{equation}\label{value:f1:f2}
 f_1=0,\quad f_2=-i2^2\xi^{-\frac{5}{4}}dv^{\frac{1}{2}}B_0^{-1}.
\end{equation}

For $\eta$ on the Stokes line $\arg\eta=0$ and $\eta\rightarrow0$,
from \eqref{relation zeta and eta zero} it immediately follows
that $\zeta^2\thicksim \frac{1}{4\eta}$, then we have $\arg\zeta\sim\pi$. Therefore,
$\arg({\rm e}^{\pi i/4}\sqrt{2\xi}\zeta)\thicksim\frac{5\pi}{4}$ and $\arg({\rm e}^{-\pi i/4}\sqrt{2\xi}\zeta)\thicksim\frac{3\pi}{4}$ for $\xi>0$.
From \eqref{Second order equation asy} we have $F^{-1/4}\thicksim 2^{\frac{3}{2}}\eta$ as $\eta\rightarrow0$ for large $\xi$. Since
$(\zeta^2-\alpha^2)^{1/4}\thicksim\zeta^{1/2}$ as $\eta\rightarrow0$, by using the appropriate asymptotic formulas of $D_{\nu}(z)$ in
(\ref{Asymptotic behavior of D z}), we obtain from \eqref{Denote of tilde phi1} and \eqref{relation zeta and eta zero} that
\begin{align}\label{asy:tilde:phi:zero1}
 \tilde{\phi}_\nu 
 \thicksim \eta\Big(C_0\eta^{-\frac{1}{4}}{\rm e}^{-\frac{i}{8\eta}\xi}
 -\frac{\sqrt{2\pi}}{\Gamma(-\nu)}{\rm e}^{i\pi\nu}D_0\eta^{\frac{1}{4}}{\rm e}^{\frac{i}{8\eta}\xi}
 \Big) ,~{\rm as}~\eta\rightarrow0,
\end{align}
where
\begin{align}
  C_0&=2^{\frac{\nu+3}{2}}\exp\big(\frac{5\pi i}{4}\nu+\frac{1}{2}\pi
  i\big)\exp\big(\frac{i}{4}\xi+\frac{\nu}{2}\ln\xi\big)(1-b)^{-\frac{1}{2}},\label{coeff:C0}\\
  D_0&=2^{1-\frac{\nu}{2}}\exp\big(-\frac{5\pi i}{4}\nu-\frac{3}{4}\pi
  i\big)\exp\big(-\frac{i}{4}\xi-\frac{\nu+1}{2}\ln\xi\big)(1-b)^{\frac{1}{2}}.\label{coeff:D0}
\end{align}
Similarly, from \eqref{Denote of tilde phi2} and \eqref{relation zeta and eta zero} it follows that
\begin{align}\label{asy:tilde:phi:zero2}
 \tilde{\phi}_{-\nu-1} 
 \thicksim \eta\Big(D_0{\rm e}^{\frac{\pi i}{2}(\nu+1)}\eta^{\frac{1}{4}}{\rm e}^{\frac{i}{8\eta}\xi}
 -\frac{\sqrt{2\pi}}{\Gamma(\nu+1)}{\rm e}^{-\frac{3\pi i}{2}\nu-\pi i}C_0\eta^{-\frac{1}{4}}{\rm e}^{-\frac{i}{8\eta}\xi}
 \Big) ,~{\rm as}~\eta\rightarrow0,
\end{align}
Substituting \eqref{asy:tilde:phi:zero1} and \eqref{asy:tilde:phi:zero2} into \eqref{asy:Phi:21:orgin} and \eqref{asy:Phi:22:orgin}, respectively,
 we obtain
\begin{align}
\Psi_{21}^{(0)}&\thicksim \xi^{\frac{1}{2}}q^{\frac{1}{2}}\Big[
\big(\delta_1-\delta_2\frac{\sqrt{2\pi}}{\Gamma(\nu+1)}{\rm e}^{-\frac{3\pi i}{2}\nu-\pi i}\big)C_0
\eta^{-\frac{1}{4}}{\rm e}^{-\frac{i}{8\eta}\xi}\nonumber\\
&\quad+\big(\delta_2{\rm e}^{\frac{\pi i}{2}(\nu+1)}
-\delta_1\frac{\sqrt{2\pi}}{\Gamma(-\nu)}{\rm e}^{\pi i\nu}\big)D_0\eta^{\frac{1}{4}}{\rm e}^{\frac{i}{8\eta}\xi}\Big],~~{\rm as}~\eta\rightarrow0,\label{asy:behavior:Psi:21:zero}\\
\Psi_{22}^{(0)}&\thicksim \xi^{\frac{1}{2}}q^{\frac{1}{2}}\Big[
\big(\delta_3-\delta_4\frac{\sqrt{2\pi}}{\Gamma(\nu+1)}{\rm e}^{-\frac{3\pi i}{2}\nu-\pi i}\big)
C_0\eta^{-\frac{1}{4}}{\rm e}^{-\frac{i}{8\eta}\xi}\nonumber\\
&\quad+\big(\delta_4{\rm e}^{\frac{\pi i}{2}(\nu+1)}
-\delta_3\frac{\sqrt{2\pi}}{\Gamma(-\nu)}{\rm e}^{\pi i\nu}\big)D_0\eta^{\frac{1}{4}}{\rm e}^{\frac{i}{8\eta}\xi}\Big],~~{\rm as}~\eta\rightarrow0.\label{asy:behavior:Psi:22:zero}
\end{align}
By using the boundary condition \eqref{Psi boundary at zero}, we have
\begin{align}
\frac{\tilde{d}}{\sqrt{h-1}}\lambda^{\frac{1}{4}}{\rm e}^{\frac{i}{8\lambda}}&\thicksim
  \xi^{\frac{1}{2}}q^{\frac{1}{2}}\Big[
\big(\delta_1-\delta_2\frac{\sqrt{2\pi}}{\Gamma(\nu+1)}{\rm e}^{-\frac{3\pi i}{2}\nu-\pi i}\big)C_0
\eta^{-\frac{1}{4}}{\rm e}^{-\frac{i}{8\eta}\xi}\nonumber\\
&\quad+\big(\delta_2{\rm e}^{\frac{\pi i}{2}(\nu+1)}
-\delta_1\frac{\sqrt{2\pi}}{\Gamma(-\nu)}{\rm e}^{\pi i\nu}\big)D_0\eta^{\frac{1}{4}}{\rm e}^{\frac{i}{8\eta}\xi}\Big],
\label{matching:zero:21}\\
 \frac{-i\sqrt{h}}{\tilde{d}\sqrt{h-1}}\lambda^{-\frac{1}{4}}{\rm e}^{-\frac{i}{8\lambda}}&\thicksim
 \xi^{\frac{1}{2}}q^{\frac{1}{2}}\Big[
\big(\delta_3-\delta_4\frac{\sqrt{2\pi}}{\Gamma(\nu+1)}{\rm e}^{-\frac{3\pi i}{2}\nu-\pi i}\big)C_0
\eta^{-\frac{1}{4}}{\rm e}^{-\frac{i}{8\eta}\xi}\nonumber\\
&\quad+\big(\delta_4{\rm e}^{\frac{\pi i}{2}(\nu+1)}
-\delta_3\frac{\sqrt{2\pi}}{\Gamma(-\nu)}{\rm e}^{\pi i\nu}\big)D_0\eta^{\frac{1}{4}}{\rm e}^{\frac{i}{8\eta}\xi}\Big].\label{matching:zero:22}
\end{align}
From \eqref{matching:zero:21} and \eqref{matching:zero:22} it follows that
\begin{equation}\label{delta2:1}
  \frac{\delta_2}{\delta_1}=\frac{{\rm e}^{\frac{3\pi i}{2}\nu+\pi i}\Gamma(\nu+1)}{\sqrt{2\pi}},~ 
  \frac{\delta_4}{\delta_3}=\frac{\sqrt{2\pi}{\rm e}^{\frac{\pi i}{2}(\nu-1)}}{\Gamma(-\nu)},~
  \delta_3=\frac{\frac{-i\sqrt{h}}{\tilde{d}\sqrt{h-1}}\xi^{-\frac{1}{4}}q^{-\frac{1}{2}}C_0^{-1}}{1
  -\frac{2\pi{\rm e}^{-\pi i\nu-\frac{3}{2}\pi i}}{\Gamma(-\nu)\Gamma(\nu+1)}}.
\end{equation}
Substituting \eqref{value:f1:f2} and \eqref{delta2:1} into \eqref{express:Q:21}, yields
\begin{align}\label{asym:Q:21}
  (Q)_{21}&=-\frac{\sqrt{2\pi}{\rm e}^{-\frac{3\pi i}{2}\nu-\pi i}}{\Gamma(\nu+1)}2^2d\tilde{d}\xi^{-1}
  (qv)^{\frac{1}{2}}B_0^{-1}C_0\frac{\sqrt{h-1}}{\sqrt{h}}\nonumber\\
  &=\frac{2^{-\frac{1}{4}}\sqrt{\pi}{\rm e}^{-\frac{\pi\beta }{4}}}{\Gamma(\frac{1}{2}-i\frac{\beta}{2})}\exp\Big(iS+\frac{i}{2}\xi-\frac{i\beta}{2}\ln\xi-\frac{i\beta}{2}\ln2
  +\frac{3\pi i}{4}\Big)(1+O(\xi^{-1})),
\end{align}
where $\beta=k^2$. This completes the proof of Lemma \ref{Q:a>pi}. \qed\\
\subsection{Proof of the lemma \ref{Q:a<pi}}
For the case of $a<\frac{1}{\pi}$, we applying the result \cite[(A.44)]{Sulei}
\begin{equation}\label{asy of Phi2}
  \Phi'(\xi,a)=1+\frac{2r^2}{\xi}\sin^2(2S)+O(\xi^{-2}),~{\rm as}~\xi\to\infty,
\end{equation}
where $S=\frac{1}{2}\Phi(\xi,a)$, and $a\in L$, here $L$ is defined as follows,
\begin{equation}\label{Defin L}
 L=\left\{a\Big|\Phi(x)=\Phi(x,a) {\rm ~is~increasing~as~} x\to\infty, \exists x>0,\Phi(x,a)>\frac{\pi}{2}\right\}.
\end{equation}
By the same argument of approximating equation \eqref{Second order equation2} in last subsection, after a carefully calculation, we obtain that
 \begin{equation}\label{Second order equation asy2}
  \frac{d^2\phi}{d\eta^2}=-\xi^2\tilde{F}(\xi,\eta)\phi,
\end{equation}
where
\begin{equation}\label{expr:2F}
  \tilde{F}(\xi,\eta)=\frac{1}{64}\big(1-\frac{1}{\eta^2}\big)^2
  +\frac{1}{\xi}\Big[-\frac{1}{4\eta^2}\big(r^2-\frac{i}{\eta}\big)
  +\frac{i}{8\eta}\big(1-\frac{1}{\eta^2}\big)\big(\frac{1}{2}-\frac{1}{b\eta-1}\big)\Big]+F_2(\eta)O(\frac{1}{\xi^2}),
\end{equation}
with $b=i\tan S$, and $\tilde{F}_2(\eta)=1+\frac{1}{\eta^2}+\frac{1}{\eta^3}$. We note that there are two coalescing turning
points near $1$, and two close to $-1$. Here we will only concern with the two turning points, say $\hat{\eta}_1$ and $\hat{\eta}_2$, near $1$. When $\hat{\eta}_j$ approach to 1, by using \eqref{Second order equation asy}, we get the asymptotic formulas for the two turning points
\begin{equation}\label{2 turning points}
  \hat{\eta}_j^{-1}=1\pm 2\xi^{-1/2}\sqrt{r^2-i}(1+o(1)),~j=1,2,
\end{equation}
 which  coalescing to $1$ when $\xi\rightarrow\infty$, and the Stokes' curves defined by $\Im(\xi(\eta+\frac{1}{\eta}))=0$. Assuming that $\xi\in\mathbb{R}^+$, then, the Stokes lines of the solution $\phi$ to \eqref{Second order equation asy2} are the positive and the negative real lines in the $\eta$ plane. Similar to \eqref{Define alpha} and \eqref{Define zeta}, if we define $\hat{\alpha}$ and  $\vartheta(\eta)$ by
\begin{equation}\label{2Define alpha}
  \frac{1}{2}\pi i\hat{\alpha}^2=\int_{-\hat{\alpha}}^{\hat{\alpha}}(\tau^2-\alpha^2)^{1/2}d\tau=\int_{\hat{\eta}_1}^{\hat{\eta}_2}F^{1/2}(\xi,s)d s,
\end{equation}
 and
\begin{equation}\label{2Define zeta}
  \int_{\alpha}^{\vartheta}(\tau^2-\alpha^2)^{1/2}d\tau=\int_{\eta_2}^\eta F^{1/2}(\xi,s)ds,
\end{equation}
respectively, then we have the following lemma which is similar to Lemma \ref{Uniform asymptotic theorem}.

\begin{lemma}\label{Uniform asymptotic lem}~
There exist constants $\hat{c}_1,~\hat{c}_2$ such that
 \begin{equation}\label{2Uniform approximate of phi}
   \phi(\eta,\xi)=\Big(\frac{\vartheta^2-\alpha^2}{\tilde{F}(\xi,\eta)}\Big)^{\frac{1}{4}}
\Big\{[\hat{c}_1+o(1)]D_{\nu}({\rm e}^{\pi i/4}\sqrt{2\xi}\vartheta)
+[\hat{c}_2+o(1)]D_{-\nu-1}({\rm e}^{-\pi i/4}\sqrt{2\xi}\vartheta)\Big\},
 \end{equation}
as $\xi\rightarrow+\infty$ uniformly for $\eta$ on the Stokes cures, where $D_{\nu}(z)$ and $D_{-\nu-1}(z)$ are solutions of the parabolic cylinder equation and $\nu$ defined by $\nu=-\frac{1}{2}+\frac{1}{2}i\xi\hat{\alpha}^2$.
\end{lemma}

Moreover, an argument similar to the one used in Lemma \ref{lem relation zeta and eta infty} obtain the asymptotic behaviors of $\vartheta(\eta)$ as $\eta\to+\infty$ and $\eta\to0$ for large $\xi$, we state those results as follows without proof.
\begin{lemma}\label{lem theta and eta infty}
For large $\xi$ and $\eta$,
\begin{equation}\label{theta and eta infty}
  \frac{1}{2}\vartheta^2=\frac{\alpha^2}{2}\ln\vartheta+\frac{1}{8}(\eta+\frac{1}{\eta})-\frac{1}{4}+\frac{i}{4\xi}\ln\eta + \frac{i}{2\xi}\ln(1-b^{-1})-
   \frac{i}{\xi}\ln2+o(\xi^{-1}),
\end{equation}
where $b=i\tan S$, and
\begin{equation}\label{asy tidalpha}
  \hat{\alpha}^2=\frac{r^2-i}{\xi}+o(\frac{1}{\xi})\quad as~\xi\to\infty.
\end{equation}
\end{lemma}
\begin{remark}
By using of the definition of $\nu$ in Lemma \ref{Uniform asymptotic lem} and \eqref{asy tidalpha}, we get the approximate value
\begin{equation}\label{2asym:nu}
  \nu=\frac{ir^2}{2}+o(1),\quad {\rm as}~\xi\to+\infty
\end{equation}
for the order of the parabolic cylinder function $D_{\nu}({\rm e}^{\pi i/4}\sqrt{2\xi}\vartheta)$ in \eqref{2Uniform approximate of phi}.
\end{remark}
\begin{lemma}\label{lem theta and eta zero}
When $\eta\to 0$, for large $\xi$, such $\xi\eta=o(1)$, then holds
\begin{equation}\label{theta and eta zero}
  \frac{1}{2}\vartheta^2=\frac{\vartheta^2}{2}\ln\vartheta+\frac{1}{8}(\eta+\frac{1}{\eta})-\frac{1}{4}-\frac{i}{4\xi}\ln\eta + \frac{i}{2\xi}\ln(1-b)-
   \frac{i}{\xi}\ln2-\frac{1}{2}\pi i\alpha^2+o(\xi^{-1}).
\end{equation}
\end{lemma}

The proofs of those lemmas are analogous to that in Lemma \ref{lem relation zeta and eta infty} and will not be include here.
Now we are in a position to prove the Lemma \ref{Q:a<pi}.

\noindent {\bf{Proof of Lemma \ref{Q:a<pi}}}.
Based on the Lemma \ref{Uniform asymptotic lem}, \ref{lem theta and eta infty} and \ref{lem theta and eta zero}, by suitable modification to the proof of Lemma \ref{Q:a>pi}, we can show that the entry $(2,1)$ of the connection martix $Q$ as $\xi\to+\infty$ for $a<\frac{1}{\pi}$ have the asymptotic behavior
\begin{align}\label{2asym:Q:21}
  (Q)_{21}
  =\frac{ir2^{-\frac{3}{4}}\sqrt{\pi}{\rm e}^{-\frac{\pi i}{2}\nu}}{\Gamma(\nu+1)}\exp\Big(-iS+i\frac{1}{2}\xi+\nu\ln\xi+\nu\ln2\Big)(1+O(\xi^{-1})).
\end{align}
Substituting \eqref{2asym:nu} into \eqref{2asym:Q:21} and denoting $\beta=r^2$, then we obtain \eqref{Q2:-infty}. This completes the proof of Lemma \ref{Q:a<pi}. \qed
\section*{Acknowledgements}
The work of Zhao-Yun Zeng was supported in part by the National Natural
Science Foundation of China under Grant Numbers 11571375 and 11571376, Guang
Dong Natural Science Foundation under Grant Number 2014A030313176 and Doctoral
Startup Fund of Jingangshan University under Grant Number JZB16001.
 \begin{appendix}
\section{Appendix A. Proof of Lemma \ref{lem relation zeta and eta infty} } \renewcommand{\thesection}{\Alph{section}}
The idea to prove Lemma \ref{lem relation zeta and eta infty} is to compute the asymptotic behavior of the integrals on the two hand sides of \eqref{Define zeta}. Straightforward  integration on the left-hand side of  \eqref{Define zeta} yields
\begin{equation}\label{Left hand integral}
  \int_{\alpha}^{\zeta}(\tau^2-\alpha^2)^{1/2}d\tau=\frac{1}{2}\{\zeta(\zeta^2-\alpha^2)^{1/2}-\alpha^2\ln(\zeta+(\zeta^2-\alpha^2)^{1/2})
+\alpha^2\ln\alpha\},
\end{equation}
Here, the cut for the integrand is again the line segment joining $-\alpha$ and $\alpha$, and again we take  the integration path  along the upper
edge of the cut.
Because we are going to calculate the higher-order part of the both two sides of (\ref{Define zeta}), we
will simply ignore the lower-order part in two hand sides, then we obtain that from (\ref{Left hand integral}) for large $\zeta$
\begin{equation}\label{Left hand integral Asymptotic}
  \frac{1}{2}\zeta^2-\frac{1}{2}\alpha^2\ln(2\zeta)-\frac{1}{4}\alpha^2+\frac{1}{2}\alpha^2\ln(\alpha)+O(\alpha^4\zeta^{-2})
=\int_{\eta_2}^\eta F^{1/2}(\xi,s)d s.
\end{equation}

To calculate the right-hand side of (\ref{Left hand integral Asymptotic}), we split right-hand side into two integrals respectively
\begin{equation}\label{Right hand integral}
 \int_{\eta_2}^\eta F^{1/2}(\xi,s)d s=
\Big(\int_{\eta_2}^{\eta^*}+\int_{\eta^*}^\eta\Big)F^{1/2}(\xi,s) ds:=I_1+I_2
\end{equation}
where
\begin{equation}\label{eta star}
  \eta^*=1+2T\xi^{-1/2},
\end{equation}
and $T$ is a large parameter to be specified more precisely later. In $I_1$ we taking the change
\begin{equation}\nonumber
  s=1+2t \xi^{-1/2},
\end{equation}
and ignore $F_2$, then $I_1$ can be evaluated for large $\xi$ as follows:
\begin{align}\label{Approximate value of I1}
I_1
=\frac{T^2}{2\xi}+\frac{k^2}{4\xi}+\frac{k^2}{2\xi}\ln(2T)-\frac{k^2}{4\xi}\ln(-k^2)+o(\xi^{-1}).
\end{align}
Taking $T=-\sqrt{-k^2}$ in $I_1$, we obtain (\ref{asy alpha}).

When $|\eta|\to\infty$, from  it follows \eqref{express:F:2} that $F_2=O(\eta^{-2})$, thus, we can ignore $F_2$ in the integral $I_2$ in (\ref{Right hand integral}),
 by the binomial expansion, then $I_2$ is given by
\begin{align}\label{Approximate value of I2}
I_2
\thickapprox\int_{\eta^*}^\eta\frac{1}{8}(1-\frac{1}{s^2})ds+\frac{4}{\xi}\int_{\eta^*}^\eta\frac{F_1(\xi,s)}{1-\frac{1}{s^2}}ds
\end{align}
with error term $o(\xi^{-1})$.

Since
\begin{equation}\label{split F1}
  \frac{4F_1(\xi,s)}{1-\frac{1}{s^2}}=\frac{k^2/2}{s-1}-\frac{k^2/2}{s+1}-\frac{i/4}{s}+\frac{ib/2}{bs-1},
\end{equation}
 then the second term in \eqref{Approximate value of I2} is equate to
 \begin{align}\label{second term in I2}
 \frac{4}{\xi}\int_{\eta^*}^\eta\frac{F_1(\xi,s)}{1-\frac{1}{s^2}}ds
  &=\frac{1}{\xi}\Big\{\frac{i}{4}\ln\eta-\frac{k^2}{2}\ln(2T\xi^{-\frac{1}{2}})+\frac{k^2}{2}\ln2
  -\frac{i}{2}\ln(1-b^{-1})\nonumber\\
  &\quad+O(\eta^{-1})+O(T\xi^{-\frac{1}{2}})\Big\}
 \end{align}
whilst the first term in $I_2$ is equate to
\begin{equation}\label{first term in I2}
\begin{split}
 \int_{\eta^*}^\eta\frac{1}{8}(1-\frac{1}{s^2})ds
    =\frac{1}{8}(\eta+\frac{1}{\eta})-\frac{1}{4}-\frac{T^2}{2\xi}+O(T^3\xi^{-\frac{3}{2}})
\end{split}
\end{equation}
Combining \eqref{Right hand integral}, \eqref{Approximate value of I1}, \eqref{Approximate value of I2},
\eqref{second term in I2} and \eqref{first term in I2} yields
\begin{align}\label{asy right hand integral}
   \int_{\eta_2}^\eta F^{1/2}(\xi,s)d s
   &=\frac{1}{8}(\eta+\frac{1}{\eta})-\frac{1}{4}+\frac{k^2}{4\xi}+\frac{i}{4\xi}\ln\eta-\frac{i}{2\xi}\ln(1-b^{-1})+\frac{k^2}{2\xi}\ln2\nonumber\\
   &\quad-\frac{k^2}{2\xi}\ln\frac{\sqrt{-k^2}}{\xi^{1/2}}+O(T^3\xi^{-\frac{3}{2}})+O(T\xi^{-2})+o(\xi^{-1})
\end{align}
and so, choosing $T<\xi^{\frac{1}{6}}$ and using \eqref{asy alpha} and \eqref{Left hand integral Asymptotic}, we obtain
\eqref{relation zeta and eta infty}, which completes the proof of Lemma \ref{lem relation zeta and eta infty}.
\hfill\qed

\section{Appendix B. Proof of Lemma \ref{lem relation zeta and eta zero}}\renewcommand{\thesection}{\Alph{section}}
The idea to prove Lemma \ref{lem relation zeta and eta zero} is to compute the asymptotic behavior of the integral on the right hand side integral in \eqref{Left hand integral Asymptotic}. When $\eta\to0$, let $\eta_*=1-2T\xi^{-\frac{1}{2}}$, where $T$ is a large parameter, and split the right hand side integral in
 \eqref{Left hand integral Asymptotic} into two parts
\begin{equation}\label{Right hand integral zero}
 \int_{\eta_2}^\eta F^{1/2}(\xi,s)d s=
\Big(\int_{\eta_2}^{\eta_*}+\int_{\eta_*}^\eta\Big)F^{1/2}(\xi,s) ds=J_1+J_2.
\end{equation}
The integral $J_1$ can be calculated by the similar manner as in computing $I_1$ \eqref{Approximate value of I1}, it follows that
\begin{equation}\label{Approximate value of J1}
  J_1=\frac{T^2}{2\xi}+\frac{k^2}{4\xi}+\frac{k^2}{2\xi}\ln(-2T)-\frac{k^2}{4\xi}\ln(-k^2)+O(\xi^{-1}T^{-2})+O(T\xi^{-2}).
\end{equation}
For the second integral $J_2$, according to the expressions of $F_1$ and $F_2$ in \eqref{express:F:1} and \eqref{express:F:2}, respectively,
 we have $F_1\thicksim\eta^{-3}$,
 $F_2(\xi,\eta)\thicksim\eta^{-3}$ as $\eta\to0$. Thus we have
\begin{align}\label{Approximate value of J2}
J_2&=\int_{\eta_*}^\eta[\frac{1}{64}(1-\frac{1}{s^2})^2+\frac{1}{\xi}F_1(\xi,s)+s^{-3}O(\xi^{-2})]^{\frac{1}{2}}ds\nonumber\\
&=\int_{\eta_*}^\eta\frac{1}{8}(1-\frac{1}{s^2})ds+\frac{4}{\xi}\int_{\eta_*}^\eta\frac{F_1(\xi,s)}{1-\frac{1}{s^2}}ds
+O(\xi^{-2}\ln\eta)+O(\xi^{-1}T^{-2})+O(\xi^{-3/2}T^{-1})\nonumber\\
&=\frac{1}{8}(\eta+\frac{1}{\eta})-\frac{1}{4}-\frac{i}{4\xi}\ln\eta-\frac{i}{2\xi}\ln(1-b)+\frac{k^2}{2\xi}\ln2
-\frac{k^2}{2\xi}\ln(2T\xi^{1/2})\nonumber\\
   &\quad-\frac{T^2}{2\xi}+O(T^3\xi^{-\frac{3}{2}})+O(T\xi^{-2})+o(\xi^{-1})
\end{align}
Combining \eqref{Right hand integral zero}, \eqref{Approximate value of J1} and \eqref{Approximate value of J2}, and choosing $T<\xi^{\frac{1}{6}}$,
 we have
\begin{align}\label{asy right hand integral zero}
   \int_{\eta_2}^\eta F^{1/2}(\xi,s)d s
   &=\frac{1}{8}(\eta+\frac{1}{\eta})-\frac{1}{4}+\frac{k^2}{4\xi}-\frac{i}{4\xi}\ln\eta-\frac{i}{2\xi}\ln(1-b)+\frac{k^2}{2\xi}\ln2\nonumber\\
   &\quad-\frac{1}{2}\pi i\alpha^2-\frac{k^2}{4\xi}\ln\frac{-k^2}{\xi}+o(\xi^{-1}).
\end{align}
Substituting this into \eqref{Left hand integral Asymptotic} yields \eqref{relation zeta and eta zero}. This completes the proof.
\hfill\qed\\
\end{appendix}

\end{document}